\documentclass[a4paper,11pt,oneside]{article}

\usepackage{type1cm}     
\usepackage{makeidx}     
\usepackage{graphicx}       
\usepackage{multicol}        
\usepackage[bottom]{footmisc}
\usepackage[numbers]{natbib}
\usepackage{authblk}
\usepackage[top=2cm, bottom=2cm, left=2cm, right=2cm]{geometry}

\usepackage{amsmath,amsthm,amsfonts,amssymb,amscd}
\usepackage{mathtools}
\makeindex           

\begin{document}

\title{Space-time error control using a partition-of-unity 
dual-weighted residual method applied to low mach number combustion}

\author[1]{Jan P. Thiele}
\author[1,2]{Thomas Wick}

\affil[1]{Leibniz University Hannover, Institut of Applied Mathematics, Welfengarten 1, 30167 Hannover, Germany}

\affil[2]{Universit\'e Paris-Saclay, ENS Paris-Saclay, LMT -- Laboratoire de M\'ecanique et Technologie, 91190 Gif-sur-Yvette, France}

\maketitle

\begin{abstract}
In this work, a space-time scheme for goal-oriented a posteriori 
error estimation is proposed. The error estimator is evaluated using 
a partition-of-unity dual-weighted residual method. As application, 
a low mach number combustion equation is considered. In some numerical 
tests, different interpolation variants are investigated, while 
observing convergence orders and effectivity indices between true 
errors (obtained on a sufficiently refined mesh) and the error estimator.\\
\textbf{Keywords:}\\
low mach number combustion; dual weighted residuals; finite elements; adaptivity\\
\textbf{AMS:}\\
35K58, 49M29, 65N30, 65N50, 90A25
\end{abstract}
\section{Introduction}
\label{thiele:sec:1}
This work is devoted to space-time goal-oriented a posteriori 
error control. Such space-time schemes for error estimation and adaptivity in time, space, or both, are of current interest
with various applications in parabolic problems \cite{SchVe08},
incompressible Navier-Stokes equations \cite{Besier2009,BeRa12},
dynamic Signorini, obstacle, and hyperbolic problems \cite{BlRaSchroe09,BlRaSchroe08,Rade09}
and fluid-structure interaction \cite{Fai17,FaiWi18}.
 
Our method is based on prior work, in which the dual-weighted error 
estimator is realized within a weak formulation using 
a partition-of-unity \cite{RiWi15}. We note that 
another weak realization is achieved with the so-called filtering 
approach \cite{BraackErn03}, which was already applied 
in \cite{SchVe08,Besier2009,BeRa12} to space-time error control and adaptivity.

However, the extension of 
the partition-of-unity (PU) localization to space-time settings 
has not yet been established in the published literature.
We provide a detailed algorithmic derivation of the error estimator 
and discuss the important ingredients. As numerical example, we consider 
a nonlinear low mach number combustion problem. A key part is the 
backward-in-time running adjoint problem. One difficulty is that, 
due to Galerkin orthogonality, the adjoint problem must contain 
higher order information (see \cite{BeRa01}) in the primal 
error part and vice vera in the adjoint error part. 
Recently, for stationary settings a new class of algorithms 
could be established by using interpolation information in a smart way
\cite{EndtLaWi20_smart}. 
For our proposed space-time setting we investigate 
the performance by simply choosing different combinations for the spatial finite elements, such as 
low-order $cG(1)/cG(1)$ and high-order equal elements $cG(2)/cG(2)$ as well
as the natural approach $cG(1)/cG(2)$ (low order primal and higher order adjoint for the 
primal error part). These choices are investigated with respect to 
their convergence properties and evaluation of the effectivity indices.
We will, however, limit ourselves to the variation of finite element orders in space and use the equal order 
approach for the temporal discretization.
We also notice that some preliminary results on space-time adaptivity with the 
PU-DWR method are published in \cite{ThiWi21}.

The outline of this work is as follows: In Section 
\ref{thiele:sec:2}, the low mach number combustion equations 
are introduced and their weak formulation is provided.
Next, in Section \ref{thiele:sec:3}, the discretization with 
finite elements is described. In the main Section \ref{thiele:sec:4}
a space-time PU-DWR error estimator is derived in detail.
Finally, in Section \ref{thiele:sec:5} some numerical experiments 
are reported, that show the performance of our developments.
The code for these simulations is based on our extension of the package \texttt{dwr-diffusion} \cite{kocher.etal2019}
to solve nonlinear problems. The package itself uses \texttt{deal.II} \cite{arndt.etal2020} as the finite element library.
We conclude our work in Section \ref{thiele:sec:6}.

\section{The Low Mach Number Combustion Equations}
\label{thiele:sec:2}
The nonlinear parabolic problem we want to investigate describes a combustion reaction under the low Mach number hypothesis.
Under that hypothesis the dimensionless temperature $\theta$ and the concentration of the combustible species $Y$ 
are not influencing the fluid velocity field. 

For the special but important case of $v=0$ all convection terms vanish and $\theta$ and $Y$ are only influenced by diffusion and 
by the reaction mechanism in which $Y$ combusts and raises the temperature.
For constant diffusion coefficients we arrive at the following set of equations
\begin{eqnarray}
 \partial_t\theta -\Delta\theta &= \omega(\theta,Y) \quad\text{in } \Omega\times(0,T),\\
 \partial_t Y     -\frac{1}{Le}\Delta Y     &= -\omega(\theta,Y) \quad\text{in }\Omega\times(0,T),
\end{eqnarray}
where the reaction is described by Arrhenius law
\begin{equation}
  \omega(\theta,Y) = \frac{\beta^2}{2Le}Ye^{\frac{\beta(\theta-1)}{1+\alpha(\theta-1)}}.
\end{equation}
The parameters are the Lewis number $Le = 1$, the gas expansion $\alpha = 0.8$ and the nondimensional activation energy $\beta = 10$.

Part of the boundary $\Gamma_R\subset\Gamma\coloneqq\partial\Omega$ will be cooled.
This can be described by a Robin condition $\partial_n\theta=-k\theta$. 
Following the standard procedure, we obtain the following space-time variational formulation of our problem.
Find $u=(\theta,Y)$ such that
\begin{eqnarray}
 A(u,\phi) = (\partial_t\theta,\phi^\theta) + (\nabla\theta,\nabla\phi^\theta) + (\partial_t Y,\phi^Y) + (\nabla Y,\nabla\phi^Y)
 \label{thiele:eq:variational}\\
 + (\omega(\theta,Y),\phi^Y-\phi^\theta) + \int\limits_0^T \int\limits_{\Gamma_R}k\theta\phi^\theta \mathrm{d}s\mathrm{d}t = 0 \eqqcolon F(\phi)
 \;\forall\phi=(\phi^\theta,\phi^Y)\nonumber
\end{eqnarray}
where $(f,g)$ describes the space-time scalar product defined as
\[
(f,g):= \int\limits_0^T\int\limits_\Omega fg\mathrm{d}x\mathrm{d}t.
\]
As the homogeneous Neumann condition on $\Gamma_N$ is a natural condition it does not appear in the variational formulation.
The inhomogeneous Dirichlet conditions on $\Gamma_D$ are imposed as usual and inserted into the finite element spaces.

\section{Discretization}
\label{thiele:sec:3}
As we want to use different finite elements for the temporal and the spatial test- and trial functions we will start by 
partitioning $[0,T]$ into $M$ subintervals $I_n = (t_{n-1},t_n]$, with the discretization parameter $k = t_n-t_{n-1}$.
In time we will use piecewise constant dicontinuous elements 
$\phi_k(t)\in dG(0)$.
To be able to use different refined meshes over time, so called dynamic meshes, we will 
discretize $\Omega$ on each subinterval by a triangulation $\mathcal{T}_h^n$.
Using quadrilaterals (in two-dimensional configurations) for the spatial triangulation, we can use continuous finite element functions $\phi_h(x)\in cG(p)$ of order $p$ 
as test functions and trial functions.
The fully discrete equations on each subinterval are then obtained by using $\phi{kh}\in dG(0)cG(p)$ as test functions and trial functions in 
\eqref{thiele:eq:variational}.
For a more detailed look at the discretization and the corresponding finite element spaces see \cite{SchVe08}.

\section{Space-Time PU-DWR Error Estimation}
\label{thiele:sec:4}
Denoting our quantity of interest by the goal functional $J(u)$, 
we obtain the Lagrange functional for minimizing the error in said quantity as
\begin{equation}
 \mathcal{L}(u,z) = J(u)+F(z)-A(u,z).
\end{equation}
As a first order optimality condition we obtain the KKT (Karush-Kuhn-Tucker) system and with it an auxiliary adjoint problem. In summary, we then have
\begin{eqnarray}
 \mathcal{L}'_u(u)(\psi,z) = J_u'(u)(\psi)-A_u'(u)(\psi,z) \stackrel{!}{=} 0\; \text{ (adjoint problem)}\\
 \mathcal{L}'_z(z)(u,\phi) = F(\phi)-A(u,\phi) \stackrel{!}{=} 0\; \text{ (primal problem)}
\end{eqnarray}
Note that for nonlinear problems $A(u,z)$ is a semilinearform that is linear in $z$ and that $F(z)$ is always a linear form.
Thus, the directional derivative in direction $\phi$ w.\,r.\,t. to $z$, i.\,e. $A'_z(z)(u,\phi)$ coincides with $A(u,\phi)$.
The same holds for linear goal functionals and linear problems resulting in the dual problem $A(\psi,z) = J(\psi)$.
Also note that the adjoint problem obtained by this derivation applies the temporal derivative to the test function.
To rectify this, a partial integration in time is applied to $(\partial_t \psi,z)$, yielding $(\psi,-\partial_t z)$. 
This results in a problem that runs backwards in time and has a final condition instead of an initial condition.
\subsection{Error Estimation}
Following Proposition (2.1) in \cite{BeRa01} we obtain the error representation
\begin{equation}
 J(u)-J(u_{kh}) = \frac12\mathcal{L}'_z(z_{kh})(u_{kh},z-z_{kh})+\frac12\mathcal{L}'_u(u_{kh})(u-u_{kh},z_{kh})+\mathcal{R},
 \label{thiele:eq:error_rep}
\end{equation}
where $\mathcal{R}$ is a higher order remainder term.
In many cases it is sufficient to approximate the error by only computing the primal residual i.e.
\begin{align}
 J(u)-J(u_{kh}) &\approx \mathcal{L}'_z(z_{kh})(u_{kh},z-z_{kh}) \nonumber \\ 
&= F(z-z_{kh})-A(u_{kh},z-z_{kh})\eqqcolon\rho_{kh}(u_{kh},z-z_{kh}),
\end{align}
which is also called primal error estimator. Subsequently the second term in the error representation is called adjoint error estimator 
$\rho_{kh}^{*}(u_{kh})(u-u_{kh},z_{kh})$.
Introducing the semidiscrete solutions $u_k$ and $z_k$ which are still continuous in space the primal error estimator can be further split 
into a temporal estimator $\rho_k$ and a spatial estimator $\rho_h$
\begin{align}
 J(u)-J(u_{kh}) &= [J(u)-J(u_k)]+[J(u_k)-J(u_{kh})]\nonumber \\
&\approx \rho_k(u_k,z-z_k) + \rho_h(u_{kh},z_k-z_{kh}).
\end{align}

\subsection{Practical Evaluation}
As the exact solutions $u$ and $z$ are unknown, we have to further approximate them to calculate the error estimators.
For the temporal primal estimator, we will construct a piecewise linear solution $i_k^{(1)}z$ on each grid point by linear interpolation 
between the piecewise constant solutions $z_n$ and $z_{n-1}$ 
in the interval $I_n$. 

For the spatial estimator we will look at three different approaches. 
The simplest approach is calculating $u_{kh}$ with $dG(0)cG(1)$ elements and $z_{kh}$ with $dG(0)cG(2)$ elements. 
Then, we assume $z_{kh}$ to be the approximation of the exact solution and interpolate it down into $cG(1)$ in space obtaining $i_h^{(2,1)}z_{kh}$. 
This interpolation should be included in most finite element packages.

For also calculating the dual estimator we also need an approximation for $u$, which can be obtained by approximating $u_{kh}$ and $z_{kh}$ with
$dG(0)cG(2)$ elements. Using the same interpolation as before we can approximate the discrete solutions as 
$\tilde{u}_{kh}\coloneqq i_h^{(2,1)}u_{kh}$ and $\tilde{z}_{kh}\coloneqq i_h^{(2,1)}z_{kh}$, while the exact solutions are approximated
by $u_{kh}$ and $z_{kh}$.

As this approach can be quite memory intensive, both $u_{kh}$ and $z_{kh}$ can be solved using $dG(0)cG(1)$ elements.
The approximation for $z$ can then be obtained by combining neighbouring $cG(1)$ elements into one large $cG(2)$ patch with the 
operator $i_{2h}^{(2)}$. The operator and the requirements for the mesh are described in \cite{BraackErn03}.
Using those interpolations we obtain the following primal estimators for the different approaches
\begin{eqnarray}
 \eta_k^{cG(1)/cG(1)} = F(i_k^{(1)}z_{kh}-z_{kh})-A(u_{kh},i_k^{(1)}z_{kh}-z_{kh}),\\
 \eta_h^{cG(1)/cG(1)} = F(i_{2h}^{(2)}z_{kh}-z_{kh})-A(u_{kh},i_{2h}^{(2)}z_{kh}-z_{kh}),
\end{eqnarray}
for $u_{kh}\in dG(0)cG(1)$ and $z_{kh}\in dG(0)cG(1)$,
\begin{align}
 \eta_k^{cG(1)/cG(2)} &= F(i_k^{(1)}i_h^{(2,1)}z_{kh}-i_h^{(2,1)}z_{kh})- A(u_{kh},i_k^{(1)}i_h^{(2,1)}z_{kh}-i_h^{(2,1)}z_{kh}),\\
 \eta_h^{cG(1)/cG(2)} &= F(z_{kh}-i_h^{(2,1)}z_{kh}) - A(u_{kh},z_{kh}-i_h^{(2,1)}z_{kh}),
\end{align}
for $u_{kh}\in dG(0)cG(1)$ and $z_{kh}\in dG(0)cG(2)$,
\begin{align}
 \eta_k^{cG(2)/cG(2)} &= F(i_k^{(1)}\tilde z_{kh}-\tilde z_{kh})- A(\tilde u_{kh},i_k^{(1)}\tilde z_{kh}-\tilde z_{kh}),\\
 \eta_h^{cG(2)/cG(2)} &= F(z_{kh}-\tilde z_{kh}) - A(\tilde u_{kh},z_{kh}-\tilde z_{kh}),
\end{align}
for $u_{kh}\in dG(0)cG(2)$ and $z_{kh}\in dG(0)cG(2)$.

The corresponding dual estimators are obtained by the same interpolation operators, but applied to the primal solution $u_{kh}$
and inserted into the adjoint problem.
\begin{eqnarray}
 \eta_k^{*cG(1)/cG(1)} = J'_u(u_{kh})(i_k^{(1)}u_{kh}-u_{kh}) - A'_u(u_{kh})(i_k^{(1)}u_{kh}-u_{kh},z_{kh})\\
 \eta_h^{*cG(1)/cG(1)} = J'_u(u_{kh})(i_{2h}^{(2)}u_{kh}-u_{kh}) - A'_u(u_{kh})(i_{2h}^{(2)}u_{kh}-u_{kh},z_{kh})
\end{eqnarray}
\begin{eqnarray}
 \eta_k^{*cG(1)/cG(2)} = J'_u(u_{kh})(i_k^{(1)}u_{kh}-u_{kh}) - A'_u(u_{kh})(i_k^{(1)}u_{kh}-u_{kh},i_{h}^{(2,1)}z_{kh})\\
 \eta_h^{*cG(1)/cG(2)} = J'_u(u_{kh})(i_{2h}^{(2)}u_{kh}-u_{kh}) - A'_u(u_{kh})(i_{2h}^{(2)}u_{kh}-u_{kh},i_{h}^{(2,1)}z_{kh})
\end{eqnarray}
\begin{align}
 \eta_k^{*cG(2)/cG(2)} &= J'_u(\tilde u_{kh})(i_k^{(1)}\tilde u_{kh}-\tilde u_{kh}) - A'_u(\tilde u_{kh})(i_k^{(1)}\tilde u_{kh}-\tilde u_{kh},\tilde z_{kh}),\\
 \eta_h^{*cG(2)/cG(2)} &= J'_u(\tilde u_{kh})(u_{kh}-\tilde u_{kh}) - A'_u(\tilde u_{kh})(u_{kh}-\tilde u_{kh},\tilde z_{kh}).
\end{align}

\subsection{Variational PU Localization}
For use in adaptive refinement we need to obtain indicators $\eta^n_K$ or $\eta^n_i$ for each cell $K$ or DoF $i$ on the time interval $I_n$, such that
\begin{equation}
 \eta = \sum\limits_{n=1}^M\sum\limits_{K\in\mathcal{T}_h^n}\eta_K^n = \sum\limits_{n=1}^M\sum\limits_{i\in \mathcal{T}_h} \eta_i^n.
\end{equation}
We propose a DoF-wise partition of unity (PU) $\chi_i^n$, with 
\begin{equation}
\sum\limits_{n=1}^M\sum\limits_{i\in \mathcal{T}_h} \chi_i^n \equiv 1, 
\end{equation}
the simplest choice is $\chi_i^n\in dG(0)cG(1)$. 
Effectively, this leads to a spatial PU per time step, that is identical to the approach of \cite{RiWi15} for stationary problems.
The estimators are obtained by multiplying the directions in the derivatives of the Lagrangian with the PU, which leads to the 
localization of the original error representation \eqref{thiele:eq:error_rep}:
\begin{eqnarray}
 2[J(u)-J(u_{kh})]_i^n \coloneqq \mathcal{L}'_z(z_{kh})(u_{kh},(z-z_{kh})\chi_i^n)+\mathcal{L}'_u(u_{kh})((u-u_{kh})\chi_i^n,z_{kh}),\\
 J(u)-J(u_{kh}) = \mathcal{R}+\sum\limits_{n=1}^M\sum\limits_{i\in \mathcal{T}_h} [J(u)-J(u_{kh})]_i^n.\;
\end{eqnarray}
Finally, inserting the PU into the estimators described in the previous subsection yields the error indicators for each space-time DoF.
\section{Numerical Example}
\label{thiele:sec:5}
\begin{figure}[h!]
\centering
\includegraphics[scale=.45]{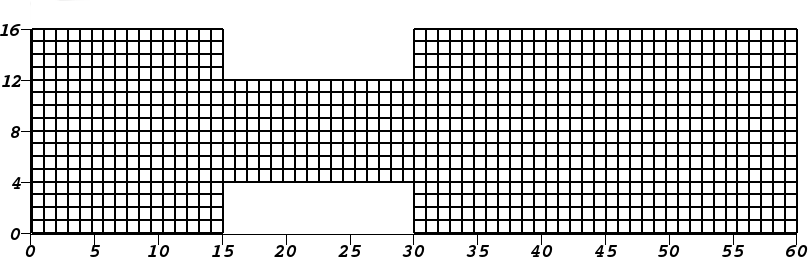}
\caption{initial grid with $N=896$ cells and $1970 (cG(1)^2)$ and $7522 (cG(2)^2)$ degrees of freedom.}
\label{thiele:fig:initial_grid}     
\end{figure}

In our numerical example, we solve 
the equations described in Section \ref{thiele:sec:2} on the geometry shown in Fig.\,\ref{thiele:fig:initial_grid}.
Here, the left edge of the domain $\Gamma_D$ is kept at a constant temperature $\theta = 1$ without any combustible species $Y=0$.
The recessed area between $x = 15$ and $x = 30$ is the cooled Robin boundary $\Gamma_R$, with $\partial_n\theta = -0.1\theta$ and $\partial_n Y = 0$. 
On the rest of the boundary $\Gamma_N$ homogeneous Neumann conditions are set. 

The initial conditions are described by
\begin{eqnarray}
 \theta_0(x,z) = \begin{cases} 1,\qquad x\leq 9, \\ e^{9-x},\quad x>9\end{cases}\\
 Y_0(x,z) = \begin{cases} 0,\qquad\qquad\quad x\leq 9, \\ 1-e^{Le(9-x)},\quad x>9\end{cases}
\end{eqnarray}
The functional of interest is the space-time averaged reaction rate
\begin{equation}
 J(u) = \frac{1}{T|\Omega|}\int\limits_0^T\int\limits_\Omega \omega(\theta,Y) \mathrm{d}x\,\mathrm{d}t,
\end{equation}
with final time $T=60$.

The initial grid is solved with $k = 0.234375$, resulting in $M=256$ time intervals.

\subsection{Comparison of Error Estimators}
To compare the estimators without influence of the adaptivity, the following simulations were done with 
global refinement in space and time. Tables \ref{thiele:primal_estimators} - \ref{thiele:full_estimators} show the results for the
error $J(u)-J(u_{kh})$ in comparison with the primal, adjoint and full estimators respectively. 
Since we use a different number of degrees of freedom for each approach, we decided to base the comparison on the number of time steps $M$ and
the number of spatial cells $N$. As a rough approximation the solution $(\theta,Y)\in cG(1)^2$ needs $2N$ and the solution $(\theta,Y)\in cG(2)^2$ needs $8N$ 
spatial degrees of freedom.
As a $cG(1)$ representation of the solution is inserted as $u_{kh}$ in \eqref{thiele:eq:error_rep}, the error $J(u)-J(u_{kh})$ is evaluated in either the 
$cG(1)/cG(1)$ or the $cG(1)/cG(2)$ case. Note that the error would be larger for the projection of the $cG(2)$ solution into $cG(1)$ for the $cG(2)/cG(2)$ approach, 
as that interpolation simply uses the values at the vertex DoFs with the respective $cG(1)$ basis functions and discards all other parts of the solution.

\begin{table}[h!]
\caption{Primal estimators for different global refinement levels.}
\label{thiele:primal_estimators}
\begin{tabular}{p{1cm}p{1cm}p{3.2cm}p{3.2cm}p{3.2cm}p{3.2cm}}
\hline\noalign{\smallskip}
$M$ & $N$ & $J(u)-J(u_{kh})$ & $\eta_{kh}^{cG(1)/cG(1)}$ & $\eta_{kh}^{cG(1)/cG(2)}$ & $\eta_{kh}^{cG(2)/cG(2)}$\\
\noalign{\smallskip}\hline\noalign{\smallskip}
$256$  & $896$   & $1.07197741e-02$ & $1.92058654e-03$ & $7.61264827e+06$ & $5.30861697e-04$\\
$512$  & $3584$  & $2.48965242e-03$ & $6.43743177e-04$ & $4.08892798e-01$ & $2.34439679e-04$\\
$1024$ & $14336$ & $5.67024544e-04$ & $2.61386699e-04$ & $1.32470408e-03$ & $2.02825765e-04$\\
$2048$ & $57344$ & $1.11743216e-04$ & $1.08270274e-04$ & $1.62070260e-04$ & $1.02680294e-04$\\
\noalign{\smallskip}\hline\noalign{\smallskip}
\end{tabular}
\end{table}

\begin{table}[h!]
\caption{Adjoint/dual estimators for different global refinement levels.}
\label{thiele:adjoint_estimators}
\begin{tabular}{p{1cm}p{1cm}p{3.2cm}p{3.2cm}p{3.2cm}p{3.2cm}}
\hline\noalign{\smallskip}
$M$ & $N$ & $J(u)-J(u_{kh})$ & $\eta_{kh}^{cG(1)/cG(1)}$ & $\eta_{kh}^{cG(1)/cG(2)}$ & $\eta_{kh}^{cG(2)/cG(2)}$\\
\noalign{\smallskip}\hline\noalign{\smallskip}
$256$  & $896$   & $1.07197741e-02$ & $1.19648379e-03$ & $1.82646766e+08$ & $1.21421097e-03$\\
$512$  & $3584$  & $2.48965242e-03$ & $8.33534627e-04$ & $2.79437497e+00$ & $7.95366447e-04$\\
$1024$ & $14336$ & $5.67024544e-04$ & $4.25107914e-04$ & $4.21358177e-03$ & $4.21309997e-04$\\
$2048$ & $57344$ & $1.11743216e-04$ & $2.74717107e-04$ & $3.83350535e-04$ & $1.70568548e-04$\\
\noalign{\smallskip}\hline\noalign{\smallskip}
\end{tabular}
\end{table}

\begin{table}[h!]
\caption{Full estimators for different global refinement levels.}
\label{thiele:full_estimators}
\begin{tabular}{p{1cm}p{1cm}p{3.2cm}p{3.2cm}p{3.2cm}p{3.2cm}}
\hline\noalign{\smallskip}
$M$ & $N$ & $J(u)-J(u_{kh})$ & $\eta_{kh}^{cG(1)/cG(1)}$ & $\eta_{kh}^{cG(1)/cG(2)}$ & $\eta_{kh}^{cG(2)/cG(2)}$\\
\noalign{\smallskip}\hline\noalign{\smallskip}
$256$  & $896$   & $1.07197741e-02$ & $1.37575686e-03$ & $9.20610356e+07$ & $8.20354657e-04$\\
$512$  & $3584$  & $2.48965242e-03$ & $5.59566023e-04$ & $1.60163388e+00$ & $5.14903063e-04$\\
$1024$ & $14336$ & $5.67024544e-04$ & $3.32992866e-04$ & $2.54117001e-03$ & $2.74961514e-04$\\
$2048$ & $57344$ & $1.11743216e-04$ & $1.77760869e-04$ & $2.55543805e-04$ & $1.25485484e-04$\\
\noalign{\smallskip}\hline\noalign{\smallskip}
\end{tabular}
\end{table}
Comparing the results over all tables, we see that the equal order approaches perform relatively similar and better than the
mixed order approach. Especially on lower refinement levels the $cG(1)/cG(2)$ results are orders of magnitude above the actual error.
On closer inspection the adjoint solutions get larger with each time step, which leads to the estimator being dominated by the indicators on the
first few time intervals. As the codes for solving the $cG(2)$ adjoint problems for $cG(2)/cG(2)$ and $cG(1)/cG(2)$ basically only differ 
in the $u_{kh}$ that is inserted in the assembly of the right hand side and the nonlinear part of the matrix, we surmise that
the errors from inserting a lower order solution get amplified with each time step. 
This would also explain why this approach does not fail for stationary problems even on coarse meshes.
In previous simulations we also saw that pairing the $cG(1)/cG(2)$ approach with solving the adjoint with $cG(1)$ elements in time 
led to worse results than the temporal equal order approach, even for the linear heat equation with the $L_2$ error as functional of interest.
For adaptivity on dynamic meshes this is of course a considerable problem, as it is advisable to start with a coarse mesh 
to only capture the local behaviour of the solution/functional at each time interval. 

When comparing the tables with each other, we see that for this problem the primal error estimator on itself performs better than the adjoint error estimator 
and is comparable to the full estimator. 
Overall, the $cG(1)/cG(1)$ approach is preferable as it is considerably cheaper to calculate compared to the $cG(2)/cG(2)$ approach for which 
multiple linear systems need to be solved with $\approx 8N$ unknowns (one solve for the adjoint and several solves for the primal Newton solver steps)
instead of $\approx 2N$ unknowns. Additionally, the primal solution vectors have to be kept either in RAM or on hard disk as the adjoint is
solved backwards in time, so the $cG(2)/cG(2)$ leads to a considerable increase in memory or storage demand. 
From a computational standpoint one can see why the $cG(1)/cG(2)$ approach would be a nice tradeoff between accuracy and memory demand
as only a single linear solve per time step has to performed on the larger set of unknowns.

\subsection{Adaptive Results}
\label{thiele:subsec:52}

\begin{figure}[h!]
\centering
\includegraphics[scale=.8]{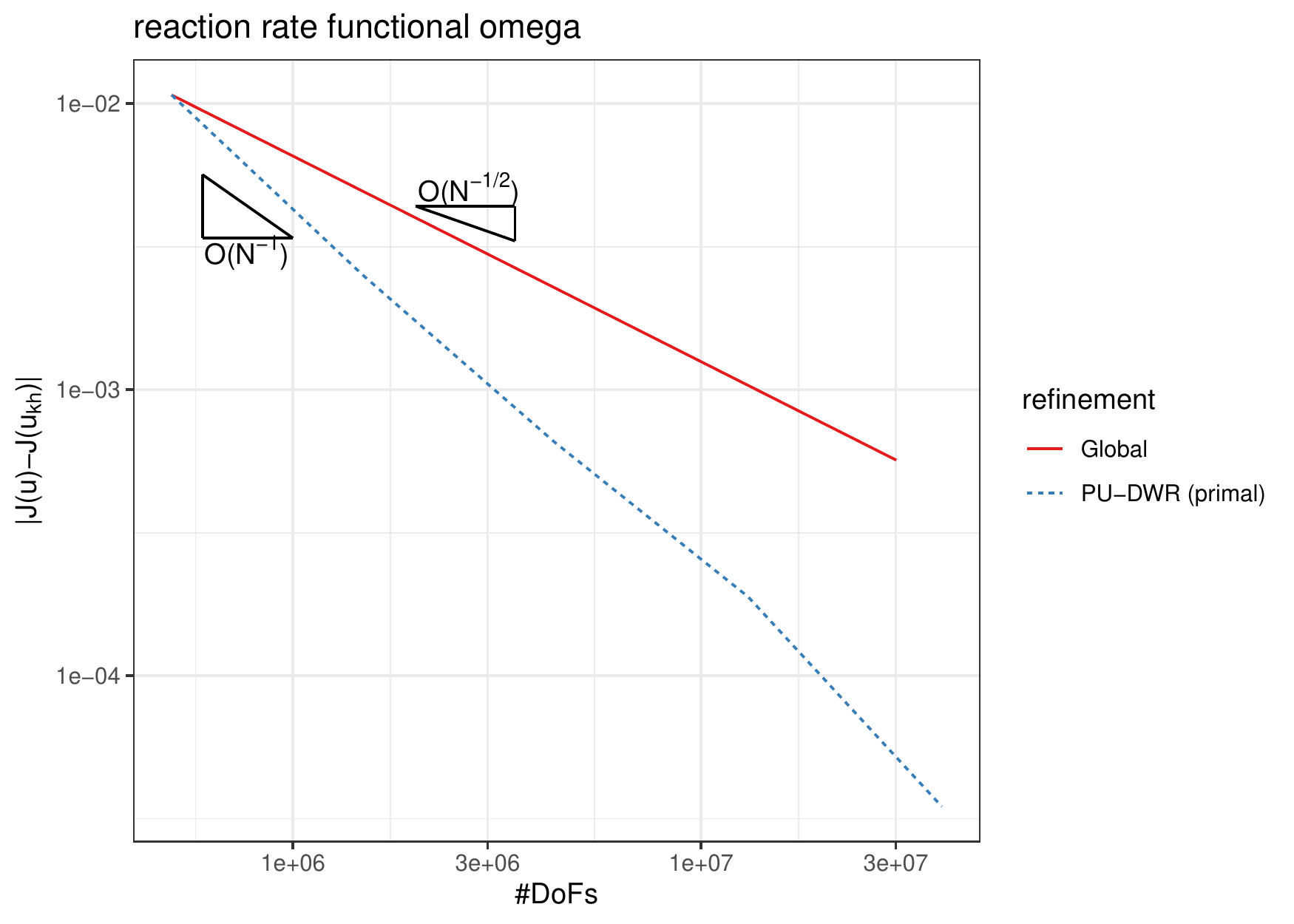}
\caption{Comparison of the actual errors for global refinement vs. adaptive refinement with the primal $cG(1)/cG(1)$ estimator with 
marking $50\%$ of the time intervals and $33\%$ of the spatial cells on each interval for refinement.}
\label{thiele:fig:error_convergenze}     
\end{figure}

When comparing the different estimators in the $cG(1)/cG(1)$ approach the primal estimator is closest to the actual error,
so we decided to use this estimator as a basis for an adaptive simulation.
As we have DoF-based indicators we compute cell-wise indicators to use build-in refinement strategies in \texttt{deal.II}.
These are obtained by simply adding the four spatial indicators of the cell vertices. 
As a refinement strategy we chose fixed fraction marking in which the indicators are sorted and the $x\%$ of cells with the 
largest indicators are marked for refinement. 
The same strategy is applied to the time intervals for which the indicators are calculated as the sum over all temporal DoF-indicators 
on the corresponding spatial triangulation.
As fractions we chose $50\%$ for the temporal and $33\%$ for the spatial refinement which leads to roughly $1.5M$ time intervals and $2N$ spatial cells
per time interval compared to $2M$ and $4N$ for global refinement.
Figure \ref{thiele:fig:error_convergenze} shows that the exact error converges faster for adaptive refinement, when comparing the 
number of primal DoFs. 

\begin{figure}[h!]
\centering
\includegraphics[scale=.24]{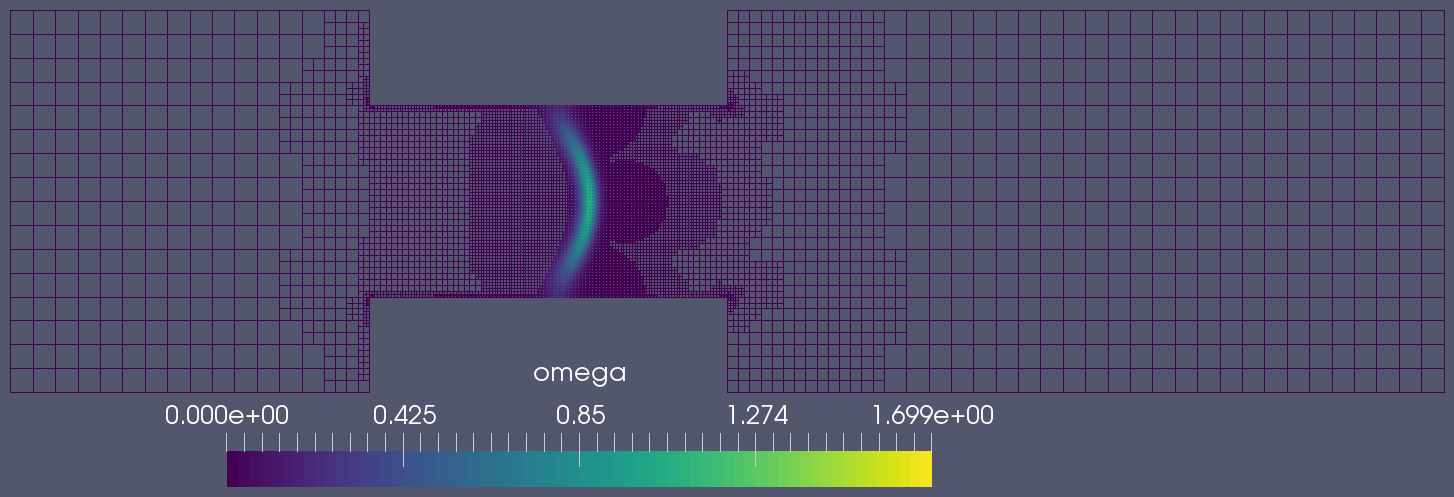}
\caption{reaction rate $\omega$ at $t=20$}
\label{thiele:fig:omega_t_20}       
\end{figure}

\begin{figure}[h!]
\centering
\includegraphics[scale=.24]{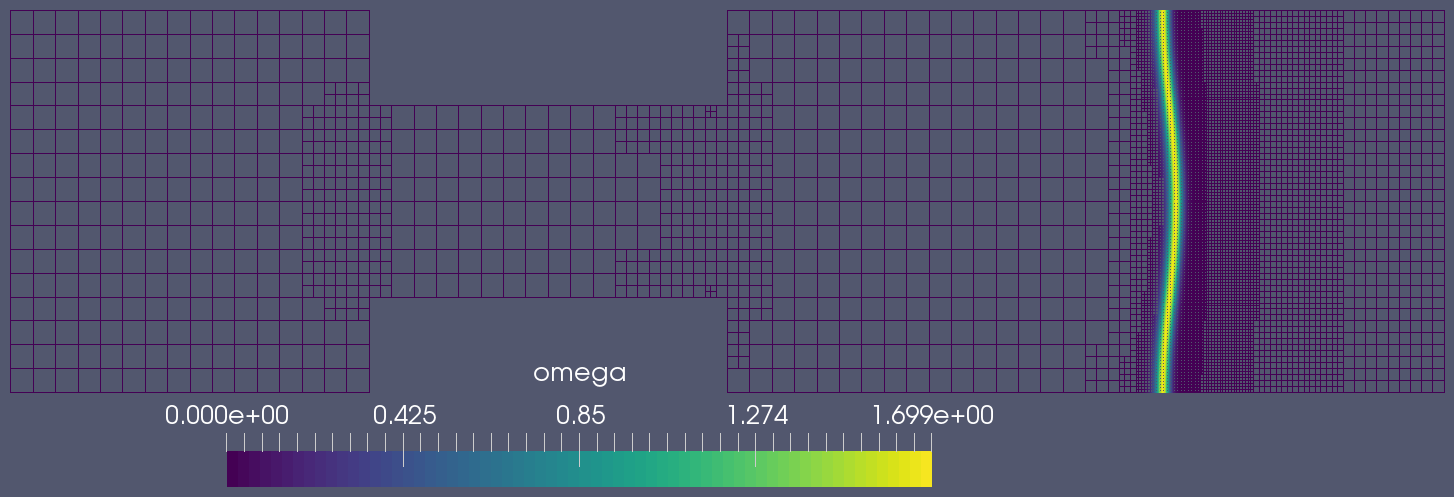}
\caption{reaction rate $\omega$ at $t=60$}
\label{thiele:fig:omega_t_60}       
\end{figure}

To see if our novel localization approach works well in capturing the local behaviour of the goal functional, Figures 
\ref{thiele:fig:omega_t_20} and \ref{thiele:fig:omega_t_60} show the evolution of the reaction rate $\omega$ 
over $[0,T]$ and the corresponding meshes. In all timesteps the combustion reaction is captured well by the fine cells.Additionally, for 
time steps after the flame passed the cooled rods, there is also some refinement around the sharp corners, which is to be expected.

\mbox{}
\newpage

\newpage
\section{Conclusions}
\label{thiele:sec:6}
In this work, we developed a space-time goal-oriented a posteriori 
error estimator using a newly developed partition-of-unity dual-weighted 
residual localization. As model problem we considered a nonlinear 
low Mach number combustion problem. Specific emphasis was on different 
space-time finite element combinations for the primal and adjoint subproblems.
Therein, we detected a better performance for equal-order combinations 
of $cG(1)/cG(1)$ and $cG(2)/cG(2)$ type in comparison to a $cG(1)/cG(2)$ finite element.
The reason has not yet been fully understood by us and needs further 
future investigations whether algorithmic or mathematical problems are 
the reason. Finally, some illustrations of locally adaptive meshes 
show that the error indicators obtained by our proposed method yield
excellent findings in terms of resolving the local flame front.

\bibliographystyle{abbrv}
\bibliography{lit.bib}
\end{document}